\documentclass[11pt, reqno, oneside]{amsart}

\usepackage{hyperref}
\usepackage{comment}
\usepackage{stmaryrd}
\usepackage{cleveref}

\usepackage[margin=1 in,letterpaper]{geometry}

\theoremstyle{definition}
\newtheorem{nul}{}[section]
\newtheorem{dfn}[nul]{Definition}
\newtheorem{rmk}[nul]{Remark}
\newtheorem{exm}[nul]{Example}
\newtheorem{cnv}[nul]{Convention}

\theoremstyle{plain}
\newtheorem{thm}[nul]{Theorem}
\newtheorem{prop}[nul]{Proposition}
\newtheorem{lem}[nul]{Lemma}
\newtheorem{cor}{Corollary}[nul]
\newtheorem*{thm*}{Theorem}
\newtheorem*{prop*}{Proposition}
\newtheorem*{cor*}{Corollary}
\newtheorem*{lem*}{Lemma}
\newtheorem*{qst*}{Question}

  \newtheoremstyle{TheoremNum}
        {\topsep}{\topsep}              
        {\itshape}                      
        {}                              
        {\bfseries}                     
        {.}                             
        { }                             
        {\thmname{#1}\thmnote{ \bfseries #3}}
  \theoremstyle{TheoremNum}

\usepackage{amsmath}

\usepackage{tikz}
\usetikzlibrary{matrix,arrows,decorations}
\usepackage{tikz-cd}

\usepackage{adjustbox}
\usepackage[utf8]{inputenc}
\usepackage[english]{babel}

\usepackage[backend=biber, style=alphabetic]{biblatex}
\begin{document}

\title{On the Bousfield Classes of $\mathrm{H}_\infty$-ring spectra}

\author{Jeremy Hahn}
\email{jhahn01@mit.edu}
\address{Department of Mathematics, Massachusetts Institute of Technology, 77 Massachusetts Ave, Cambridge, MA 02139}

\begin{abstract}
We prove that any $K(n)$-acyclic, $D_p$-ring spectrum is $K(n+1)$-acyclic, affirming an old conjecture of Mark Hovey. 
\end{abstract}

\maketitle
\tableofcontents


\section{Introduction}

A bedrock result of chromatic homotopy theory is that any $K(h)$-acyclic, $p$-local finite spectrum is $K(h-1)$-acyclic.  Our goal here is to prove that $\mathbb{E}_\infty$-ring spectra enjoy the opposite phenomenon:
\begin{thm} \label{thm:intro-main}
Suppose that an $\mathbb{E}_{\infty}$-ring spectrum $R$ is $K(h)$-acyclic at some prime $p$. Then $R$ is also $K(h+1)$-acyclic at $p$.
\end{thm}
In fact, our arguments require much less than a full $\mathbb{E}_{\infty}$-structure:

\begin{dfn}
For a fixed prime number $p$, a $D_p$-algebra $R$ is a spectrum $R$ equipped with:
\begin{enumerate}
\item A unit map $\mathbb{S} \to R$ and a multiplication map $R \otimes R \to R$, making $R$ into a homotopy commutative and associative ring spectrum.
\item A factorization of the $p$-fold multiplication map $R^{\otimes p} \to R$ through the projection $R^{\otimes p} \to (R^{\otimes p})_{hC_p}$, such that the diagram
\[
\begin{tikzcd}[column sep=huge]
(\mathbb{S}^{\otimes p})_{hC_p} \arrow{d} \arrow{r}{(1^{\otimes p})_{hC_p}} & (R^{\otimes p})_{hC_p} \arrow{d} \\
\mathbb{S} \arrow{r}{1} & R
\end{tikzcd}
\]
commutes up to homotopy.  Here, $1$ is the unit map, and the left-hand vertical map is part of the natural $\mathbb{E}_{\infty}$-ring structure on the sphere $\mathbb{S}$.
\end{enumerate}
\end{dfn}

Any $\mathbb{E}_{\infty}$-ring $R$ is naturally a $D_p$-algebra, and any $\mathrm{H}_{\infty}$-ring $R$ admits a $D_p$-algebra structure.  In particular, \Cref{thm:intro-main} is a corollary of the following stronger result:



\begin{thm} \label{thm:intro-main-ii}
Suppose that a spectrum $R$ admits a $D_p$-algebra structure.  For any height $h \ge 0$ at the prime $p$, if $R$ is $K(h)$-acyclic then $R$ is also $K(h+1)$-acyclic.
\end{thm}

\begin{rmk}
If a spectrum $R$ admits a homotopy unital multiplication, then $R$ is $K(h)$-acyclic if and only if it is $T(h)$-acyclic \cite[Lemma 2.3]{LMMT}, so \Cref{thm:intro-main-ii} can also be read as a statement about telescopic localization.
\end{rmk}

\begin{rmk}
The $h=0$ case of \Cref{thm:intro-main-ii} was proved by Mathew--Naumann--Noel \cite[Theorem 2.1]{MNN}\footnote{While Mathew--Naumann--Noel state their theorem for $\mathrm{H}_{\infty}$-ring spectra, their proof uses only $D_p$-algebra structure.}--it is known as the May nilpotence conjecture. Together with Clausen, these authors found spectacular applications of the May nilpotence conjecture to descent questions in algebraic $K$-theory \cite{Descent}.
\end{rmk}

In order to prove \Cref{thm:intro-main-ii}, it suffices to consider not $D_p$-algebras in spectra, but rather $K(h+1)$-local $D_p$-$E$-algebras, where $E$ is the height $h+1$ Morava $E$-theory with $$\pi_*E \cong \mathbb{Z}_{p^{h+1}}\llbracket u_1,u_2,\cdots,u_{h} \rrbracket[u^{\pm}].$$
We will study a $K(h+1)$-local $D_p$-$E$-algebra $R$ by means of the weight $p$ \emph{power operation} 
\[P:\pi_0(R) \to R^0(BC_p),\] which takes $x:S^0 \to R$ to the homotopy class of the composite
\[BC_p \simeq (\mathbb{S}^{\otimes p})_{hC_p} \xrightarrow{(x^{\otimes p})_{hC_p}} (R^{\otimes p})_{hC_p} \to R.\]
The power operation $P$ is multiplicative, but not additive. 
Nonetheless, the composition of $P$ with the quotient map $R^0(BC_p) \to R^0(BC_p)/(\mathrm{tr})$ is additive.  Here, $(\mathrm{tr})$ denotes the \emph{transfer ideal}, which is cut out by a certain element $\frac{[p](z)}{z} \in R^0(BC_p)$ that is related to the $p$-series of the formal group law on $\pi_0(E)$.
In particular, for each $x \in \pi_0(R)$ the image of $P(x)$ in the quotient ring $R^0(BC_p)/(p,\mathrm{tr})$ depends only on $x$ modulo $p$.
One of our main technical results, which may be of substantial independent interest, is a higher chromatic analog of the preceding sentence.  For each $1 \le k \le h$, we determine a quotient of $R^{BC_p}$ in which $P(x)$ depends only on the value of $x$ modulo $(p,u_1,\cdots,u_k)$:

\begin{thm} \label{thm:poweropintro}
Suppose that $R$ is a $K(h+1)$-local $D_p$-$E$-algebra and $k \ge 1$. Then, for any $x \in \pi_0(R)$, the image of $P(x)$ in $\pi_0\left(R^{BC_p}/\left(p,u_1,\cdots,u_k,\frac{[p](z)}{z^{p^{k+1}}}\right) \right)$ depends only on the image of $x$ in $\pi_0(R/(p,u_1,\cdots,u_k))$.
\end{thm}

For a more detailed version of this statement, see \Cref{sec:ModPowOp}, where $\frac{[p](z)}{z^{p^{k+1}}}$ is denoted $g_{k+1}(z)$.  Note that $[p](z)$ is divisible by $z^{p^{k+1}}$ only after modding out by $(p,u_1,\cdots,u_k)$.

\begin{rmk}
In private communication, Nathaniel Stapleton has asked whether there exist power operations for the cohomology theory of spaces given by the spectrum $E/(p,u_1,\cdots,u_{k})$.  In other words, if $X$ is a space, Stapleton asks whether there is a natural map
\[\left(E/(p,u_1,\cdots,u_{k})\right)^{0}(X) \to \left(E^{BC_p}/\left(p,u_1,\cdots,u_{k},\frac{[p](z)}{z^{p^{k+1}}}\right)\right)^{0}(X).\]
The above result shows that such a power operation can be defined on any $x \in \left(E/(p,u_1,\cdots,u_{k})\right)^{0}(X)$ that arises as the image of a class in $E^0(X)$.
\end{rmk}

Let us briefly describe how we prove \Cref{thm:intro-main-ii}. Suppose that $R$ is a $K(h+1)$-local $D_p$-$E$-algebra.  If $R$ is $K(h)$-acyclic, then for some value of $n$ it must be the case that $(u_{h})^n \in \pi_0(E)$ maps to zero in $\pi_0(R/(p,u_1,\cdots,u_{h-1}))$.  Using \Cref{thm:poweropintro}, we check that if $(u_h)^n$ maps to zero, then so does $(u_{h})^{n-1}$.  Iterating this argument, we eventually learn that one maps to zero in $\pi_0(R/(p,u_1,\cdots,u_{h-1}))$, which implies (since $R$ is $K(h+1)$-local) that $R$ is the zero ring.  Such a proof is very much analogous to the one of the May nilpotence conjecture in \cite{MNN}.

\begin{rmk}
In the long time since the first version of this paper was posted, Clausen--Mathew--Naumann--Noel \cite{CMNN} and Land--Mathew--Meier--Tamme \cite{LMMT} found applications of \Cref{thm:intro-main} to foundational descent and purity results in algebraic $K$-theory. As part of their arguments, Clausen--Mathew--Naumann--Noel make the elegant observation that, for any $K(h+1)$-local $\mathbb{E}_{\infty}$-ring $R$, $R$ is $K(h)$-acyclic if and only if $R^{tC_p}$ is $K(h)$-acyclic, where $R^{tC_p}$ denotes the Tate construction for the trivial $C_p$ action.  We record as \Cref{thm:CMNNvariant} that this is also true of $D_p$-algebras.  Thus, if a $K(h+1)$-local $D_p$-algebra $R$ has $L_{K(h)}R^{tC_p}=0$, then $R=0$.
\end{rmk}

\begin{rmk}
In their forthcoming work on the chromatic Nullstellensatz, Burklund--Schlank--Yuan prove that any $\mathbb{E}_{\infty}$-ring is either $K(h+1)$-acyclic or maps to some height $h+1$ Lubin--Tate theory; this in particular strengthens \Cref{thm:intro-main}. In combination with work of Allen Yuan \cite{YuanRed}, the Nullstellensatz implies chromatic redshift occurs in the algebraic $K$-theory of any $\mathbb{E}_{\infty}$-ring.
\end{rmk}

\begin{rmk}
A version of \Cref{thm:intro-main-ii} was first conjectured by Mark Hovey, where it appears as Miscellaneous Problem $2$ in his 1999 list of unsolved problems in algebraic topology \cite{HoveyProblems}.
\end{rmk}

\subsection*{Acknowledgements}
I heartily thank Akhil Mathew for introducing me to this problem and pointing out its appearance on Mark Hovey's webpage.  Thanks are due to Peter May, Denis Nardin, Eric Peterson, Andrew Senger, Nathaniel Stapleton, Dylan Wilson, Allen Yuan, and especially my advisor Mike Hopkins for helpful conversations.  The author was supported by the NSF Graduate Fellowship under Grant DGE-1144152.

The original version of this paper was the first preprint I ever produced, and I hope the subsequent six years of experience have improved both my technical and expository powers. In particular, I thank Robert Burklund, Akhil Mathew, and Andrew Senger for bringing to my attention a mistake in Section $2$ of the original version of this preprint.

\section{Morava $E$ theory and its weight $p$ power operation}

We fix throughout the remainder of this paper a prime $p$ and an integer $h>0$, and by default all spectra will be implicitly $p$-localized.  By convention, we let $E=E_{h+1}$ denote the height $h+1$ Lubin--Tate theory associated to the Honda formal group over $\mathbb{F}_{p^{h+1}}$.  We call this $E_{h+1}$ \emph{Morava $E$-theory}, and it has homotopy groups
\[\pi_*E \cong \mathbb{Z}_{p^{h+1}} \llbracket u_1,u_2,\cdots,u_h \rrbracket [u^{\pm}],\]
where $|u_i|=0$ and $|u|=-2$.  By convention, we sometimes write $u_0=p$.  Any other Lubin--Tate theory would work just as well for our purposes, but we use this one for the sake of concreteness.  The Goerss--Hopkins--Miller theorem equips $E$ with a canonical $\mathbb{E}_\infty$-ring structure \cite{GoerssHopkins}.

\begin{cnv}
We fix a $p$-typical complex orientation of $E$
\[BP \to E,\]
and name our generators $u_i \in \pi_0E$ such that indecomposable generators $v_i \in \pi_{2p^i-2} BP$ are carried to $u_i u^{-p^i+1}$ \cite[p. 7]{DevinatzHopkins}.  Note that we do not in any way insist that this orientation be compatible with the $E_{\infty}$-ring structure on Morava $E$-theory.  Our fixed $p$-typical complex orientation is given by a class $t \in E^2(\mathbb{CP}^{\infty})$, and it will be convenient to let $z \in E^0(\mathbb{CP}^{\infty})$ denote the product of $t$ with $u^{-1} \in \pi_{2}E$. We can then speak of the $p$-series $[p](z)$ as a class in $E^0(\mathbb{CP}^{\infty}) \cong E_0\llbracket z \rrbracket$.
\end{cnv}

\begin{rmk}
The fibration 
\[S^1 \to \mathrm{BC}_p \to \mathbb{CP}^{\infty}\]
allows one to write $E^0(BC_p)$ as a quotient of $E^0(\mathbb{CP}^{\infty}) \cong E_0\llbracket z \rrbracket$ \cite[Lemma 5.7]{HKR}.  Specifically, one has
\[E^0(BC_p) \cong E_0\llbracket z \rrbracket/[p](z).\]
\end{rmk}

\begin{dfn}
The central object of study of this paper is \emph{the weight $p$ power operation on Morava $E$-theory}, which is a multiplicative (but not additive) map
\[P:E^0 \to E^0(BC_p).\]
By definition, this map takes a class $x:S^0 \to E$ to the composite class 
\[\Sigma^{\infty}_+ BC_p = (S^0)^{\otimes p}_{hC_p} \xrightarrow{(x^{\otimes p})_{hC_p}} (E^{\otimes p})_{hC_p} \to E,\]
where the last map arises from the $D_p$-algebra structure on Morava $E$-theory.
\end{dfn}

\begin{dfn} While $P$ is not additive, and hence not a ring map, the composite
\[E^0 \stackrel{P}{\longrightarrow} E^0(BC_p) \longrightarrow z^{-1} E^0(BC_p) = \pi_0E^{tC_p}\]
is a ring map, which we will denote by $\varphi$.  In fact, $\varphi$ is $\pi_0$ of the Nikolaus--Scholze Frobenius \cite[\textsection{IV.1}]{NikolausScholze}, which is an $\mathbb{E}_\infty$ ring map
$E \to E^{tC_p}.$
By taking $\pi_*$ of the Nikolaus--Scholze Frobenius, we extend the domain of $\varphi$ to include elements of $\pi_*E$ that are not of degree $0$; the output of $\varphi$ will in general be a class in $z^{-1}E^*(BC_p)$.
\end{dfn}

\begin{lem} \label{TateFrobeniusCongruence}
For each $0 \le k \le h$, the following congruence holds in $z^{-1}E^0(BC_p)$:
\[\varphi(u_k) \equiv \frac{\varphi\left(u^{p^k-1}\right) u_k}{u^{p^k-1}} \text{ modulo } p,\cdots,u_{k-1}\]
\end{lem}

\begin{proof}
Recall that we have fixed a homotopy ring map $f:\mathrm{BP} \to E$, which on homotopy groups takes $v_i$ to $u_i u^{-p^i+1}$.
The lemma is thus implied by the statement that $\varphi(v_i) \equiv v_i$ modulo $(p,v_1,\cdots,v_{i-1})$ in $z^{-1}E^0(BC_p)$.  To see this, consider the homotopy ring map given by the composite
\[\mathrm{BP} \otimes \mathrm{BP} \xrightarrow{f \otimes f} E \otimes E \xrightarrow{can \otimes \varphi} E^{tC_p},\]
where $\mathrm{can}$ denotes the canonical map that exists because we are considering $E$ with trivial $C_p$-action. Taking homotopy groups, $v_i$ will be the image of $\eta_L(v_i) \in \pi_*(\mathrm{BP} \otimes \mathrm{BP})$, while $\varphi(v_i)$ will be the image of $\eta_R(v_i)$.  The result then follows from the fact that $\eta_L(v_i) \equiv \eta_R(v_i)$ modulo $(p,v_1,\cdots,v_{i-1})$ in $\pi_*(\mathrm{BP} \otimes \mathrm{BP})$.
\end{proof}

\begin{cor} \label{corTateFrobeniusCongruence}
For each $0 \le k \le h$, $\varphi(u_k) \equiv 0$ modulo $p,\cdots,u_k$ in $z^{-1}E^0(BC_p)$.
\end{cor}

It will be useful to formulate the above corollary in terms of $P$, instead of $\varphi$.  To do so, recall (by, e.g., the formula at the top of \cite[p.788]{GHMR}) that  the $p$-series $[p](z) \in E_0\llbracket z \rrbracket$ satisfies the equation
\[[p](z) \equiv u_kz^{p^k} \text{ modulo } p,\cdots, u_{k-1}, z^{p^k+1},\]
for each $0 \le k \le h$. We may thus make the following definition:

\begin{dfn}
For each $0<k\le h$, let $g_k(z) \in E_0\llbracket z \rrbracket$ denote the power series such that $z^{p^k} g_k(z)$ is obtained from the $p$-series $[p](z)$ by setting $p,u_1,\cdots,u_{k-1}=0$.
\end{dfn}

\begin{rmk}
By the Weierstrass preparation theorem \cite[Lemma 5.1]{HKR}, $z$ is not a zero-divisor in $E_0\llbracket z \rrbracket/(p,u_1,\cdots,u_{k-1},g_k(z))$. 
\end{rmk}

\begin{prop} \label{prop:finalsec2}
In $E^0(BC_p)$, the following congruence holds for each $0 \le k \le h-1$:
\[P(u_{k}) \equiv 0 \text{ }\mathrm{ modulo }\text{ } p,u_1,\cdots,u_{k},g_{k+1}(z).\]
\end{prop}

\begin{proof}
By \Cref{corTateFrobeniusCongruence} there exists some positive integer $i$ such that $z^iP(u_{k})$ is trivial in $E^0(BC_p)$ modulo $p,u_1,\cdots,u_{k}$.  The result now follows from the fact that $z$ is not a zero-divisor in $E^0(BC_p)/(p,u_1,\cdots,u_{k},g_{k+1}(z))$.
\end{proof}
 
\section{The weight $p$ power operation in a $D_p$-$E$-algebra} \label{sec:ModPowOp}

In the preceding section, we studied the weight $p$ power operation $P:E^0 \to E^0(BC_p)$, where $E$ is the height $h+1$ Morava $E$-theory.  In this section, we study the analogous operation for any $K(h+1)$-local $D_p$-$E$-algebra $R$:

\begin{dfn}
A $K(h+1)$-local $D_p$-$E$-algebra is a $K(h+1)$-local $E$-module $R$ equipped with the data of:
\begin{itemize}
\item Multiplication and unit maps $R^{\otimes_{\tiny{E}} 2} \to R$ and $E \to R$, making $R$ into a commutative and associative algebra object in the homotopy category of $E$-modules.
\item A factorization of the $p$-fold multiplication map $R^{\otimes_{E} p} \to R$ through the projection $R^{\otimes_E p} \to (R^{\otimes_E p})_{hC_p}$.
\end{itemize}
\end{dfn}


\begin{exm}
If $A$ is a $D_p$-algebra in spectra, then $L_{K(h+1)}\left(E \otimes A\right)$ is naturally a $K(h+1)$-local $D_p$-$E$-algebra.
\end{exm}

If $R$ is a $K(h+1)$-local $D_p$-$E$-algebra, then $R^{BC_p}$ is in particular an $E^{BC_p}$-module, so it makes sense to speak of the quotient of $R^{BC_p}$ by any sequence of elements in $E^0(BC_p)$.  With this in mind, the main result of this section is the following theorem:

\begin{thm} \label{modPowerExists}
Let $R$ be a $K(h+1)$-local $D_p$-$E$-algebra, $k \ge 1$, and suppose that $x \in \pi_0R$ has the property that $x$ maps to $0$ in $\pi_0 \left(R/(p,u_1,\cdots,u_k)\right)$. Then $P(x) \in \pi_0(R^{BC_p})$ maps to zero in $\pi_0 \left(R^{BC_p}/\left(p,u_1,\cdots,u_k,g_{k+1}(z)\right)\right)$.
\end{thm}

Let us explain for a moment the complications involved with proving \Cref{modPowerExists}.  For each integer $0 \le i \le k$, the quotient $E/u_i$ admits the structure of a homotopy $E$-algebra. The above theorem would be straightforward, even before coning off $g_{k+1}(z)$, if each $E/u_i$ furthermore admitted a $D_p$-$E$-algebra structure. Unfortunately, the $E/u_i$ fail to admit such structure: the first obstruction is the fact that $P(u_i) \in E^0(BC_p)$ does not map to $0$ in $E^0(BC_p)/u_i$.  However, $P(u_i)$ does project to $0$ in $E^0(BC_p)/(p,u_1,\cdots,u_i,g_{i+1}(z))$, and it is this fact that will allow us to prove \Cref{modPowerExists}.

The remainder of this section will be entirely devoted to the proof of \Cref{modPowerExists}.  Before continuing, it will be helpful to introduce some language from equivariant homotopy theory.  Our use of equivariant language is entirely confined to the remainder of this section.

\begin{dfn}
The category of $E$-modules with $C_p$-action, also known as the category of $C_p$-equivariant $E$-modules, is the category of functors from the groupoid $BC_p$ to the category of $E$-modules.
If $M$ is an $E$-module, we will often write $M^{\otimes_{E} p}$ to denote the $C_p$-equivariant $E$-module $M \otimes_E M \otimes_E \cdots \otimes_E M$ where the $C_p$ action permutes $p$ tensor factors.  For readability, when it is clear from context and in this section only, we will abbreviate $M^{\otimes_E p}$ by $M^{\otimes p}$.
\end{dfn}

\begin{rmk}
The symmetric monoidal category of $K(h+1)$-local $E$-modules with $C_p$-action is equivalent to the category of $K(h+1)$-local modules over the $\mathbb{E}_{\infty}$-ring spectrum $E^{BC_p}$ \cite[Corollary 5.4.4]{Ambi}.  Given a class $x \in E^0(BC_p)$, and an $E^{BC_p}$-module $M$, we may form the $E^{BC_p}$-module $M/x = E^{BC_p}/x \otimes_{E^{BC_p}} M$.  In particular, if $M$ is a $K(h+1)$-local $E$-module with $C_p$-action, and $x \in E^0(BC_p)$, it makes sense to speak of the $C_p$-equivariant $E$-module $M/x$.  
\end{rmk}

\begin{rmk}
Expanding on the above remark, given any sequence of elements $x_1,x_2,\cdots,x_m$ in $E^0(BC_p)$, and a $K(h+1)$-local $C_p$-equivariant $E$-module $M$, we may form the $C_p$-equivariant $E$-module $M/(x_1,\cdots,x_m)$.  In fact, since $E^*(BC_p)$ is concentrated in even degrees, there exists an $\mathbb{E}_1$-$E^{BC_p}$-algebra structure on $E^{BC_p}/(x_1,\cdots,x_m)$, and a choice of such $\mathbb{E}_1$-algebra structure allows us to view $C_p$-equivariant $E$-module $M/(x_1,\cdots,x_m)$ as a module over $E^{BC_p}/(x_1,\cdots,x_m)$ \cite{HahnWilsonQuotient}.
\end{rmk}

The following lemma is the key technical fact powering our proof of \Cref{modPowerExists}.  It concerns the $C_p$-equivariant $p$th tensor power of the non-equivariant $u_i$-Bockstein map $\Sigma^{-1}E/u_i \to E$.

\begin{lem} \label{key-technical-lemma}
For any $0 \le i \le h-1$, the $C_p$-equivariant $E$-module map
\[\left(\Sigma^{-1} E/u_i\right)^{\otimes p} \to E^{\otimes p} = E\]
admits a $C_p$-equivariant section after modding out by $(p,u_1,\cdots,u_i,g_{i+1}(z))$. In other words, there is a section of the $C_p$-equivariant $E/(p,u_1,\cdots,u_i,g_{i+1}(z))$-module map
\[\left(\Sigma^{-1} E/u_i\right)^{\otimes p}/(p,u_1,\cdots,u_i,g_{i+1}(z)) \to E^{\otimes p}/(p,u_1,\cdots,u_i,g_{i+1}(z)) = E/(p,u_1,\cdots,u_i,g_{i+1}(z)).\]
\end{lem}

Before proving the lemma, we recall the following useful fact:

\begin{rmk} \label{rmk-sparsity}
The $\mathbb{E}_{\infty}$-ring spectrum $E$ has an $\mathbb{F}_p^{\times}$ action preserving $\pi_0E$.  The fixed points of this action form an $\mathbb{E}_{\infty}$-ring spectrum that we denote by $\hat{E}$, which has homotopy groups concentrated in degrees that are multiples of $2p-2$. The inclusion of homotopy fixed points is an $\mathbb{E}_{\infty}$-ring map $\hat{E} \to E$, which is an isomorphism on $\pi_0$.  We use the existence of $\hat{E}$ only to prove the lemma below; the key useful property is the triviality of homotopy groups in degrees that are not multiples of $2p-2$.
\end{rmk}

\begin{proof}[Proof of \Cref{key-technical-lemma}]
To understand this point, it is (at least for the author) helpful to think of the $C_p$-equivariant $E$-module $\left(E/u_i\right)^{\otimes p}$ as a Thom spectrum.
Indeed, the non-equivariant $E$-module $E/u_i$ is the Thom spectrum of the map
\[S^1 \xrightarrow{1+u_i} \mathrm{BGL}_1(E).\]
The $\mathbb{E}_\infty$-ring structure on $E$ induces a $C_p$-equivariant map $\mathrm{BGL}_1(E)^{\times p} \to \mathrm{BGL}_1(E)$, and we can use this to make a $C_p$-equivariant map
\[(S^1)^{\times p} \xrightarrow{(1+u_i)^{\times p}} (\mathrm{BGL}_1(E))^{\times p} \to \mathrm{BGL}_1(E)\]
that has Thom spectrum exactly the $C_p$-equivariant $E$-module $\left(E/u_i\right)^{\otimes p}$.  Our goal will be accomplished upon proving that the composite map
\[(S^1)^{\times p} \to \mathrm{BGL}_1(E) \to \mathrm{BGL_1}(E/(p,u_1,\cdots,u_i,g_{i+1}(z)))\]
is $C_p$-equivariantly nullhomotopic.  In fact, at primes $p>2$ it is easier to prove the slightly stronger statement that the composite
\[(S^1)^{\times p} \to \mathrm{BGL}_1(\hat{E}) \to \mathrm{BGL_1}(\hat{E}/(p,u_1,\cdots,u_i,g_{i+1}(z)))\]
is $C_p$-equivariantly nullhomotopic.

We accomplish this by examining the natural equivariant cell decomposition of the pointed $C_p$ space $(S^1)^{\times p}$.  Since the non-equivariant space $S^1$ admits a cell decomposition with one $0$-cell and one $1$-cell, $(S^1)^{\times p}$ admits an equivariant cell decomposition with:
\begin{itemize}
\item One $0$-cell.
\item For each $0<k<p$, $\frac{\binom{p}{k}}{p}$ induced cells with boundaries $(C_p)_+ \wedge S^{k-1}$.
\item A $\rho$-cell (i.e., a cell with boundary $S^{\rho-1}$), where $\rho$ is the real regular representation of $C_p$. 
\end{itemize}

Noting that $S^{\rho}$ is built from $S^1$ by attaching induced cells of dimension at least $2$ and at most $p$, we see that it suffices to check the following three facts:
\begin{enumerate}
\item The composite 
\[S^1 \xrightarrow{\Delta} (S^1)^{\times p} \xrightarrow{(1+u_i)^{\times p}} \mathrm{BGL}_1(\hat{E}) \to \mathrm{BGL_1}(\hat{E}/(p,u_1,\cdots,u_i,g_{i+1}(z)))\] is nullhomotopic.
\item The composite 
\[(C_p)_+ \wedge S^1 \to (S^1)^{\times p} \to \mathrm{BGL}_1(\hat{E}) \xrightarrow{(1+u_i)^{\times p}} \mathrm{BGL_1}(\hat{E}/(p,u_1,\cdots,u_i,g_{i+1}(z)))\] 
is nullhomotopic.
\item For $1<k \le p$, all maps 
\[(C_p)_+ \wedge S^{k} \to \mathrm{BGL_1}(\hat{E}/(p,u_1,\cdots,u_i,g_{i+1}(z)))\] are nullhomotopic.
\end{enumerate}

Statement $(1)$ is the claim that $P(1+u_i)=1$ in $E^{0}(BC_p)/(p,u_1,\cdots,u_{i},g_{i+1}(z))$.  Since we have killed $g_{i+1}(z)$, we have killed the generator $[p](z)/z$ of the transfer ideal in $E^0(BC_p)$, and so we may write $P(1+u_i)=P(1)+P(u_i)=1+P(u_i)$.  Now we conclude $(1)$ from \Cref{prop:finalsec2}.

To see statements $(2)$ and $(3)$, we recall that a map 
\[(C_p)_+ \wedge S^k \to \mathrm{BGL_1}(\hat{E}/(p,u_1,\cdots,u_i,g_{i+1}(z)))\]
is the data of a non-equivariant map from $S^k$ to the non-equivariant ring underlying $\hat{E}/(p,u_1,\cdots,u_{i},g_{i+1}(z))$.  Because $g_{i+1}(z)$ is $u_{i+1}$ plus a multiple of $z$, the non-equivariant ring underlying $\hat{E}/(p,u_1,\cdots,u_i,g_{i+1}(z))$ is the non-equivariant quotient $\hat{E}/(p,u_1,\cdots,u_i,u_{i+1})$.  Since $(p,u_1,\cdots,u_{i+1})$ is a regular sequence in $\pi_0 \hat{E}=\pi_0 E$, we may conclude $(3)$ by the sparsity highlighted in \Cref{rmk-sparsity}.
To see statement $(2)$, we need to check that $1+u_i=1$ in $E_0/(p,u_1,\cdots,u_{i},u_{i+1})$, which follows from the fact that we have coned off $u_i$.
\end{proof}

The next corollary studies the $C_p$-equivariant $p$th tensor power of the non-equivariant iterated Bockstein
\[\Sigma^{-k} E / (p,u_1,\cdots,u_k) \to \cdots \to \Sigma^{-2} E/ (p,u_1) \to \Sigma^{-1} E / p \to E.\]
Note that this iterated Bockstein may also be described as the tensor product of the Bocksteins
\[\Sigma^{-1} E/u_i \to E\]
as $i$ ranges from $0$ to $k$.

\begin{cor} \label{cor:bigsectionexists}
For each $0 \le k \le h-1$, there is a $C_p$-equivariant section of the map
\[(\Sigma^{-k} E/(p,\cdots,u_{k}))^{\otimes p} \to E^{\otimes p} = E\]
after modding out by $(p,u_1,\cdots,u_{k},g_{k+1}(z))$.
\end{cor}

\begin{proof}
We tensor together the sections provided by \Cref{key-technical-lemma} for $i=k,k-1,\cdots,0$.  This gives a section of the map
\[(\Sigma^{-k} E / (p,\cdots,u_k))^{\otimes p} \to E\]
after first coning off by $(p,u_1,\cdots,u_{k},g_{k+1}(z))$, then coning off $(p,u_1,\cdots,u_{k-1},g_{k}(z)),$ etc.   Each of $p,u_1,\cdots,u_{k-1},g_1(z),g_2(z),\cdots, g_k(z)$ is trivial modulo $(p,u_1,\cdots,u_{k},g_{k+1}(z))$.  So the above coning off process yields a direct sum of copies of $E/(p,u_1,\cdots,u_{k},g_{k+1}(z))$, and we may project onto a single copy.
\end{proof}

\begin{proof}[Proof of \Cref{modPowerExists}]
By assumption, we are given a $D_p$-$E$-algebra $R$ and a class $x \in \pi_0(R)$ such that $x$ maps to $0$ in $\pi_0(R/(p,\cdots,u_{k}))$.
Consider now the following commutative diagram of $C_p$-equivariant $E$-modules
\[
\begin{tikzcd}
(\Sigma^{-k} E/(p,\cdots,u_k))^{\otimes p} \arrow{r}{f} & E \arrow{r}{x^{\otimes p}} & R^{\otimes p} \\
\Sigma^{-1}(\Sigma^{-k} E/(p,\cdots,u_k))^{\otimes p}/g_{k+1}(z) \arrow{u}{g} \arrow{r}{h} & \Sigma^{-1} E / g_{k+1}(z) \arrow{u}{j}
\end{tikzcd}
\]
Here, the map $j$ is the $g_{k+1}(z)$ Bockstein $\Sigma^{-1} E/g_{k+1}(z) \to E$, where $E$ is considered as a $C_p$-equivariant $E$-module with trivial action.  The map $g$ is obtained by tensoring $j$ with the $C_p$-equivariant $E$-module $(\Sigma^{-k} E/(p,\cdots,u_k))^{\otimes p}$.  The map $f$ is the $p$th tensor power of the Bockstein $\Sigma^{-k} E/(p,\cdots,u_k) \to E$.  Finally, the map $h$ is obtained by tensoring $f$ with $\Sigma^{-1} E/g_{k+1}(z)$.

The composite $x^{\otimes p} \circ f$ is trivial, because it is the $p$th tensor power of the non-equivariant composite
\[\Sigma^{-k} E/(p,\cdots,u_k) \to E \to R\]
that is trivial by assumption.  In particular, this implies that $x^{\otimes p} \circ j \circ h$ is trivial. By \Cref{cor:bigsectionexists}, $x^{\otimes p} \circ j$ becomes trivial after further coning off $p,u_1,\cdots,u_k$.

Finally, we consider the $C_p$-equivariant map $R^{\otimes p} \to R$ into $R$ with trivial action, which is given by the assumed $D_p$-$E$-algebra structure on $R$.  From the above discussion, we learn in particular that the composite
$$\Sigma^{-1}E/g_{k+1}(z) \to E=(E^{\otimes p}) \xrightarrow{x^{\otimes p}} (R^{\otimes p}) \to R \to R/(p,u_1,\cdots,u_k)$$
is nullhomotopic. If we compose all but the initial map and final maps in this chain, we obtain $P(x) \in R^0(BC_p)$.  The full composite precisely records the image of $P(x)$ in $\pi_0 (R^{BC_p}/p,\cdots,u_k,g_{k+1}(z))$.
\end{proof}

\section{The weight $p$ power operation modulo $(p,u_1,\cdots,u_{h-1})$}

In this final section, we study in greater detail the mod $(p,u_1,\cdots,u_{h-1})$ weight $p$ power operation on a $K(h+1)$-local $D_p$-$E$-algebra $R$, beginning with the special case $R=E$.  In particular, \Cref{modPowerExists} ensures that the following definition is sensible.

\begin{dfn}
Let 
\[\overline{P}:\mathbb{F}_p \llbracket u_h \rrbracket \to \mathbb{F}_p\llbracket u_h,z\rrbracket/(g_{h}(z))\]
be the unique ring homomorphism fitting into the following diagram:
\[
\begin{tikzcd}
E_0 \arrow{d} \arrow{r}{P} & E^0(BC_p) \arrow{d}\\
E_0/(p,u_1,\cdots,u_{h-1}) \arrow{r}{\overline{P}} &E^0(BC_p)/(p,u_1,\cdots,u_{h-1},g_h(z)),
\end{tikzcd}
\]
where the vertical maps are the natural quotient homomorphisms.
\end{dfn}

The codomain of $\overline{P}$ can be understood fairly explicitly, using the Weierstrass preparation theorem as in the following proposition:
\begin{prop} \label{g-exists}
There exists a polynomial $g(z) \in \mathbb{F}_p\llbracket u_h\rrbracket[z]$ such that:
\begin{itemize}
\item In the ring $\mathbb{F}_p\llbracket u_h,z\rrbracket$, $g(z)=U g_h(z)$ where $U$ is a unit.
\item $g(z)$ is monic, of degree $p^{h+1}-p^h$, and $g(z) \equiv z^{p^{h+1}-p^{h}}$ modulo $u_h$.
\item The constant term of $g(z)$ is divisible by $u_h$ but not $u_h^{2}$.
\end{itemize}
\end{prop}

\begin{proof}
Recall that $[p](z) \equiv u_hz^{p^h}$ modulo $(p,u_1,\cdots,u_{h-1},z^{p^h+1})$, and that $[p](z)$ is a unit multiple of $z^{p^{h+1}}$ modulo $(p,u_1,\cdots,u_{h})$.
The result then follows from the Weierstrass preparation theorem \cite[Lemma 5.1]{HKR}.
\end{proof}

By the first of the above bullet points, we may describe the codomain $E^0(BC_p)/(p,u_1,\cdots,u_{h-1},g_h(z))$ of $\overline{P}$ equally well as $E^0(BC_p)/(p,u_1,\cdots,u_{h-1},g(z))$. We then have the following proposition:

\begin{prop}
The ring $E^0(BC_p) / (p,u_1,\cdots,u_{h-1},g_h(z)) \cong \mathbb{F}_p\llbracket u_h \rrbracket[z]/(g(z))$ is a discrete valuation ring, with uniformizer $z$.
\end{prop}

\begin{proof}
The ring $\mathbb{F}_p\llbracket u_n \rrbracket[z]$ is a UFD, so Eisenstein's criterion applies to show that $g(z)$ is irreducible. It follows that the quotient $\mathbb{F}_p\llbracket u_n \rrbracket[z]/g(z)$ is a local domain.  Furthermore, $\mathbb{F}_p\llbracket u_n \rrbracket[z]/(z,g(z)) \cong \mathbb{F}_p$, because the constant term of $g(z)$ is divisible by $u_h$ but not $u_h^2$.  It follows that $z$ is a uniformizer.
\end{proof}

\begin{cnv}
It will be convenient to scale the discrete valuations on $\mathbb{F}_p\llbracket u_h \rrbracket$ and $\mathbb{F}_p\llbracket u_h,z \rrbracket/(g(z))$ so that $u_h$ has valuation $1$ in both rings.  This means that $z$ has valuation $\frac{1}{p^{h+1}-p^h}$, and the valuation of any non-zero element in $\mathbb{F}_p\llbracket u_h,z \rrbracket/(g(z))$ is a multiple of $\frac{1}{p^{h+1}-p^{h}}$. 
\end{cnv}

\begin{prop} \label{prop:key}
The class $\overline{P}(u_h) \in \mathbb{F}_p\llbracket u_h,z \rrbracket / g(z)$ has valuation $\frac{p-1}{p^{h+1}-p^h}$.
\end{prop}

\begin{proof}
The valuation on $\mathbb{F}_p\llbracket u_h,z \rrbracket / g(z)$ extends to a valuation on $z^{-1}\mathbb{F}_p\llbracket u_h,z \rrbracket / g(z)$.  In this latter ring, we have the equation 
\[\varphi(u_h)=u_h\frac{\varphi(u^{p^h-1})}{u^{p^h-1}}=u_h \left(\frac{\varphi(u)}{u}\right)^{p^{h}-1},\]
by \Cref{TateFrobeniusCongruence}.
It therefore suffices to prove that the valuation of $\varphi(u)/u$ is $\frac{-(p-1)}{p^{h+1}-p^h}$, since this will imply that the valuation of $\varphi(u_h)$ is
$$1-\frac{(p-1)(p^h-1)}{p^{h+1}-p^h} = \frac{p-1}{p^{h+1}-p^h}.$$

The valuation of $\frac{\varphi(u)}{u}$ does not depend on the choice of degree $-2$ unit $u$, so we may as well study a particularly convenient choice of degree $-2$ unit.  By a theorem of Ando \cite{Ando} generalized by Zhu \cite[Corollary 8.17]{ZhuNorm}, there exists an $\mathrm{H}_{\infty}$ ring homomorphism $\mathrm{MUP} \to E$, where $\mathrm{MUP}$ is the periodic complex bordism obtained as the Thom spectrum of the $J$-homomorphism $\mathrm{BU} \times \mathbb{Z} \to \mathrm{pic}(\mathbb{S})$.  Assuming that $u$ is the image of the degree $-2$ generator of $\pi_*(\mathrm{MUP})$ under such an $\mathrm{H}_{\infty}$-ring map, we have (as in \cite[pg. 339]{NikolausScholze}) that $\varphi(u)$ is $u^p$ divided by the Euler class of the reduced real regular representation of $C_p$. In symbols, recalling that $z \in E^0(BC_p)$ denotes $u^{-1}$ times the Euler class of the standard representation of $C_p$ on $\mathbb{C}$, we have that
\[\varphi(u)=\frac{u}{z([2](z))\cdots([p-1](z))}.\]
For each integer $k$ between $1$ and $p-1$, $[k](z)=kz+\mathcal{O}(z^2)$ has valuation $\frac{1}{p^{h+1}-p^h}$, and it follows that $\varphi(u)/u$ has valuation $\frac{-(p-1)}{p^{h+1}-p^h}$.
\end{proof}

\begin{cor} \label{cor:valdecrease}
Let $x \in \mathbb{F}_p \llbracket u_h \rrbracket$ denote any non-zero element of positive valuation.  Then $\overline{P}(x) \in \mathbb{F}_p \llbracket u_h,z\rrbracket/g_h(z)$ has strictly smaller valuation.
\end{cor}

\begin{proof}
Since we may write an arbitrary $x \in \mathbb{F}_p \llbracket u_h \rrbracket$ as a power of $u_h$ times a unit, it suffices to check this for powers of $u_h$, where it follows from \Cref{prop:key} since we have assumed $h>0$.
\end{proof}

\begin{cor}
Suppose that $R$ is a $K(h+1)$-local $D_p$-$E$-algebra, and that $x \in E_0$ is a non-zero element such that the image of $x$ in $\pi_0(R/(p,\cdots,u_{h-1}))$ is trivial.  Then, if $x$ has positive valuation in $E_0/(p,\cdots,u_{h-1}) \cong \mathbb{F}_p \llbracket u_h \rrbracket$, there exists an element $y \in E_0$ which:
\begin{enumerate}
 \item Also maps to $0$ in $\pi_0(R/(p,\cdots,u_{h-1}))$.
 \item Has valuation in $E_0/(p,\cdots,u_{h-1}) \cong \mathbb{F}_p \llbracket u_h \rrbracket$ strictly smaller than that of $x$.
\end{enumerate}
\end{cor}

\begin{proof}
Let $\alpha$ denote the image of $P(x)$ in $E^0(BC_p)/(p,u_1,\cdots,u_{h-1},g_h(z))$.  By \Cref{modPowerExists}, the image of $\alpha$ in $R^0(BC_p)/(p,u_1,\cdots,u_{h-1},g_h(z))$ is trivial.

By the isomorphism 
\[E^0(BC_p)/(p,u_1,\cdots,u_{h-1},g_h(z)) \cong \mathbb{F}_p\llbracket u_h \rrbracket[z]/(g(z)),\]
where $g(z)$ is monic of degree $p^{h+1}-p^{h}$, we may write 
$$\alpha=a_0+a_1z+a_2z^2+\cdots+a_{p^{h+1}-p^{h}-1} z^{p^{h+1}-p^{h}-1}$$
for some collection of coefficients $a_i \in \mathbb{F}_p\llbracket u_h \rrbracket$.

By \Cref{cor:valdecrease}, the valuation of $\alpha$ is strictly less than the valuation of $x$ in $\mathbb{F}_p\llbracket u_h \rrbracket \cong E_0/(p,\cdots,u_{h-1})$. It follows that there exists some $k$, with $0 \le k < p^{h+1}-p^{h}$, such that $a_k$ has valuation strictly lower than that of $x$.

Noting that $E^0(BC_p)/(p,u_1,\cdots,u_{h-1},g_h(z))$ is a free $E_0/(p,u_1,\cdots,u_{h-1})$ module, with basis dual to the functions that pick out the coefficients $a_i$, we may cap with the class that picks out $a_k$ to learn that $a_k$ maps to $0$ in $\pi_0(R/(p,\cdots,u_{h-1}))$.  We may then take $y$ to be any lift of $a_k \in E_0/(p,\cdots,u_{h-1})$ to a class in $E_0$.
\end{proof}

\begin{cor} \label{cor-produce-trivial-ring}
Suppose that $R$ is a $K(h+1)$-local $D_p$-$E$-algebra, and that $x \in E_0$ maps to zero in $\pi_0(R/(p,\cdots,u_{h-1}))$.  Then, if $x$ maps to a non-zero element of $E_0/(p,\cdots,u_{h-1}) \cong \mathbb{F}_p \llbracket u_h \rrbracket$, $R$ is trivial.
\end{cor}

\begin{proof}
We learn from the above that $1=0$ in the homotopy groups of the ring spectrum $R/(p,\cdots,u_{h-1})$.  Since $R$ is $K(h+1)$-local, this implies $R \simeq 0$.
\end{proof}

We now turn to the proof of the main theorem of this article, which we restate for the reader's convenience:

\begin{thm}
Suppose that $R$ is a $D_p$-algebra in spectra.  For any height $h$ at the prime $p$, if $R$ is $K(h)$-acyclic then $R$ is also $K(h+1)$-acyclic.
\end{thm}

\begin{proof}
To check that $R$ is $K(h+1)$-acyclic, it suffices to check that $L_{K(h+1)}(R \otimes E)$ is trivial.  Furthermore, since $K(h) \otimes L_{K(h+1)}(R \otimes E)$ is a module over the homotopy ring $K(h) \otimes R$, $L_{K(h+1)}(R \otimes E)$ will be $K(h)$-acyclic if $R$ is $K(h)$-acyclic.

We may therefore assume without loss of generality that $R$ is a $K(h+1)$-local $D_p$-$E$-algebra. If $R$ is $K(h)$-acyclic, then it is $T(h)$-acyclic for some telescope $T(h)=v_{h}^{-1}\mathbb{S}/(p^i_0,v_1^{i_1},\cdots,v_{h-1}^{i_{h-1}})$ \cite[Lemma 2.3]{LMMT}, so there must be some positive integer $k$ such that $(u_h)^{k} \in \pi_0E$ maps to zero in $\pi_0\left( R/(p,\cdots,u_{h-1})\right)$.  In particular, we may apply \Cref{cor-produce-trivial-ring} and learn that $R \simeq 0$.
\end{proof}

Following the arguments of \cite{CMNN}, we also have the following variant of the main theorem:
\begin{cor} \label{thm:CMNNvariant}
Suppose that $R$ is a $D_p$-algebra in spectra, and let $R^{tC_p}$ denote the Tate construction of $R$ with trivial $C_p$ action.  For any height $h$ at the prime $p$, if $R^{tC_p}$ is $K(h)$-acyclic then $R$ is also $K(h+1)$-acyclic.
\end{cor}

\begin{proof}
Note that $\left(L_{K(h+1)}(E \otimes R)\right)^{tC_p}$ is a module over $R^{tC_p}$, so the former is $K(h)$-acyclic whenever the latter is.  The result therefore follows by the combination of the proposition below with \Cref{thm:intro-main-ii}.
\end{proof}

\begin{prop} \label{thm:CMNNvariant}
Suppose that $R$ is a $K(h+1)$-local homotopy $E$-algebra, and let $R^{tC_p}$ denote the Tate construction of $R$ with trivial $C_p$ action.  If $R^{tC_p}$ is $K(h)$-acyclic, then $R$ is $K(h)$-acyclic.
\end{prop}

\begin{proof} 
We follow the arguments of Clausen--Mathew--Naumann--Noel from \cite[\textsection{4}]{CMNN}.  First, one checks that $E^{hC_p}$ is a free and finitely generated $E$-module, which follows as in many of our above arguments by the Weierstrass preparation theorem.  It follows that the natural comparison map $E^{hC_p} \otimes_E M \to M^{hC_p}$ is an equivalence for any $K(h+1)$-local $E$-module $M$, considered with trivial $C_p$-action.  This then implies that $E^{tC_p} \otimes_E M \simeq M^{tC_p}$.

It then remains to check that $K(h) \otimes_{E} E^{tC_p} \otimes_E R \simeq 0$ if and only if $K(h) \otimes_E R \simeq 0$, which follows from the explicit description
$\pi_*(K(h) \otimes_E E^{tC_p}) \cong \mathbb{F}_p \llbracket u_h \rrbracket[z^{\pm 1}][u_h^{-1}]/g(z).$
\end{proof}

\printbibliography

\end{document}